%% file: resubmission2.tex
\documentclass[twocolumn]{autart}    

\usepackage{graphicx}          
\usepackage{natbib}

\input{myIncludes}

\begin{document}

\allowdisplaybreaks

\begin{frontmatter}

\title{Counterexamples in synchronization: pathologies of consensus seeking gradient descent flows on surfaces\thanksref{footnoteinfo}} 

\thanks[footnoteinfo]{This work was funded by the UL project OptBioSys.}



\author[Paestum]{Johan Markdahl}\ead{markdahl@kth.se},    

\address[Paestum]{Luxembourg Centre for Systems Biomedicine, University of Luxembourg}  

\begin{keyword}                           
Multi-agent systems, asymptotic stability, consensus, differential geometry, optimization               
\end{keyword}                             

\begin{abstract}                          
Certain consensus seeking multi-agent systems can be formulated as gradient descent flows of a disagreement function. We study how known pathologies of gradient descent flows in Euclidean spaces carry over to consensus seeking systems that evolve on nonlinear manifolds. In particular, we show that the norms of agent states can diverge to infinity, but this will not happen if the manifold is the boundary of a convex set. Moreover, the system can be initialized arbitrarily close to consensus without converging to it, but this will not happen if the manifold is analytic. For analytic manifolds, consensus is asymptotically stable. This last result summarizes a number of previous findings in the literature on generalizations of the well-known Kuramoto model to high-dimensional manifolds.
\end{abstract}

\end{frontmatter}

\section{Introduction}

\noindent Consensus seeking systems (CSS) on nonlinear spaces appear in various applications of control theory and physics. This includes rigid-body attitude sync on $\mathsf{SO}(3)$ \citep{sarlette2009autonomous} and quantum bit sync on the Bloch sphere \citep{lohe2010quantum}. Most result on the convergence of such systems concern particular manifolds, \eg the circle \citep{scardovi2007synchronization}, the sphere  \citep{olfati2006swarms,lohe2010quantum}, the Stiefel manifolds \citep{markdahl2020high} including $\mathsf{SO}(n)$ \citep{deville2018synchronization}, and the unitary group \citep{lohe2010quantum}. A typical result is local stability of the synchronization manifold. However, as local convergence to sync is shown to hold on multiple manifolds, it begs the question of when it fails to hold. A theory of CSS on Riemannian manifolds can give a meaningful answer to this question. The literature contains some results that characterize CSS on Riemannian manifolds in terms of the known properties of the manifolds, \eg{} giving conditions for local convergence to sync \citep{tron2013riemannian}, global convergence \citep{sarlette2009consensus}, or determining when a CSS is multistable \citep{markdahl2021synchronization}. The contribution of this paper is to provide two additional results in the same spirit: (i) a CSS on a surface that is the boundary of a convex set cannot diverge to infinity and (ii) the sync manifold of a CSS on a closed analytic Riemannian manifold is asymptotically stable. We also give counterexamples to show that the assumptions we make on the manifolds are necessary for the results (i) and (ii) to hold.

\subsection{Literature review}



Synchronization on general manifolds is studied in \citep{tron2013riemannian, aydogdu2017opinion,markdahl2021synchronization}. \citet{tron2013riemannian} establish local convergence to consensus for discrete-time systems on manifolds with bounded curvature. A key aspect of their approach is the use of intrinsic geometry, \ie of not relying on the embedding the manifold in some Euclidean space. Our results differ from those of \citep{tron2013riemannian} since we consider continuous time system dynamics in an extrinsic setting. The work  \citep{aydogdu2017opinion} considers both intrinsic and extrinsic dynamics in continuous time, but do not present any convergence results. Recently, the work \citep{markdahl2021synchronization} establishes local convergence to a class of equilibria apart from synchronization. By contrast, this paper focuses on synchronization, or consensus, which we take to have the same meaning. 

For specific manifolds, there is a plethora of papers that establishes convergence to  consensus. The CSS in this paper is a generalization of the Kuramoto model on its homogeneous form with identical natural frequency parameters. Other results on generalized Kuramoto models include a literature on quantum synchronization on the Bloch sphere and generalizations thereof to $\mathsf{SO}(n)$ \citep{deville2018synchronization} and $\mathsf{U}(n)$  \citep{lohe2010quantum}. We mention just a few of the works on Kuramoto models over manifolds: the sphere \citep{crnkic2018swarms}, the ellipsoid \citep{zhu2014high}, the hyperboloid \citep{hyperboloid}, and tensors \citep{tensors}. Similar results also appear in applications involving \eg opinion consensus \citep{aydogdu2017opinion}, bio-inspired models of source-seeking and learning \citep{al2018gradient,crnkic2018swarms}, and computer science applications like time-series clustering \citep{crnkic2019data}. Some of these results on local convergence to the consensus manifold \citep{lohe2010quantum,aydogdu2017opinion,deville2018synchronization,markdahl2018tac,hyperboloid} are summarized by our asymptotical stability result for CSS on analytic manifolds.

\section{Preliminaries}

\subsection{Hypersurfaces}


A \emph{hypersurface} $\mathcal{H}^n$ is an $n$-dimensional manifold embedded in an $n+1$-dimensional Euclidean space. We assume that $\mathcal{H}^n$ can be implicitly characterized as
\begin{align}\label{eq:hypersurface}
\mathcal{H}^n=\{\ve{y}\in\R^{n+1}\,|\,c(\ve{y})=0\},
\end{align}
where $c:\mathbb{R}^n\rightarrow\R$ is a $C^1$ function. This form of the hypersurface also yields a normal from the Gauss map $\ve{n}:\mathcal{H}^n\rightarrow\mathcal{S}^n:\ve{y}\mapsto\nabla c(\ve{y})/\|\nabla c(\ve{y})\|$ where we assume that $\nabla c(\ve{y})\neq0$ for all $\ve{y}\in\mathcal{H}^n$. It can be showed that $\mathcal{H}^n$ is a manifold using the implicit function theorem.



\subsection{Gradient descent flows}

\noindent Gradient descent flows are the continuous time equivalents of gradient descent algorithms for minizing a real valued function. Let $\M$ be embedded in an Euclidean space $\R^n$. Given a potential function $V:\M\rightarrow\R$, the \emph{gradient descent flow} of $V$ on $\M$ is the dynamical system
\begin{align*}
\vd{x}=-\grad_{\ve{x}} V(\ve{x})=-\mathbf{P}_{\ve{x}}(\nabla_{\ve{x}} V(\ve{x}))
\end{align*}
where $\grad$ is the intrinsic gradient on the tangent space, $\mathbf{P}_{\ve{x}}:\R^n\rightarrow\ts[
\M]{\ve{}x}$ is an orthogonal projection on the tangent space, and $\nabla_{\ve{x}}$ is the gradient in $\R^n$ \citep{absil2009optimization}.

\begin{proposition}[Picard-Lindel\"{o}f]\label{prop:picard}
Suppose that $V:\R\rightarrow\M$ is $C^2$ on $\M$, then the flow $\vd{x}=-\grad V(\ve{x})$ has a unique solution $\ve{x}(t)$ which exists for all $t\in\R$.
\end{proposition}

\begin{proposition}\label{prop:diverge}
The potential function $V$ decreases along the solutions of $\vd{x}=-\grad V(\ve{x})$. For any integral curve, as $t\rightarrow\pm\infty$, either $\grad V(\ve{x})\rightarrow0$ or $|V(\ve{x})|\rightarrow\infty$.
\end{proposition}

For proofs, see \citet{jost2008riemannian}. A stronger result exists for gradient flows of analytic functions on analytic manifolds. The state $\ve{x}$ of such flows either diverges to infinity or converges to a singleton set as $t\rightarrow\pm\infty$ \citep{lageman2007convergence}.


\subsection{Consensus seeking systems}

\noindent Consider a multi-agent system with $N$ agents. We use a graph $\mathcal{G}=(\mathcal{V},\mathcal{E})$ to model interactions between agents. Each node $i\in\V$ corresponds to an agent and each edge $\{i,j\}\in\mathcal{E}$ corresponds to a pair of communicating agents. Items associated with agent $i$ carry the subindex $i$; we denote the state of agent $i$ by $\ma[i]{x}\in\M$, the normal of $\M$ at $\ve[i]{x}$ by $\ve[i]{n}$ \etc We call $\ve{x}:=(\ve[i]{x})_{i=1}^N\in\M^N$ a configuration of agents. For future reference we also define the cycle graph
\begin{align}\label{eq:cycle}
\mathcal{C}_N:=(\{1,2,\ldots,N\},\{\{1,2\},\{2,3\},\ldots,\{N,1\}\}).
\end{align}

The goal of a CSS is for the agents to asymptotically approach the \emph{consensus manifold}
\begin{align}
\mathcal{C}:=\{(\ve[i]{x})_{i=1}^N\in\M^N\,|\,\ve[i]{x}=\ve[j]{x},\,\forall \{i,j\}\in\E\}.\label{eq:C}
\end{align}
The set $\mathcal{C}$ is a manifold $\mathcal{C}\cong\M$ by the diffeomorphism $\mathcal{C}\rightarrow\M:(\ve[i]{x})_{i=1}^N\mapsto\ve[1]{x}$. Define the distance between a point and a set as $d(\ve{y},\mathcal{S}):=\inf_{\ve{s}\in\mathcal{S}}\|\ve{y}-\ve{s}\|_2$. A CSS strives to achieve $d((\ve[i]{x})_{i=1}^N,\mathcal{C})\rightarrow0$ as $t\rightarrow\infty$. As another measure of the distance to consensus, consider the \emph{disagreement function} $V:\M^N\rightarrow\R$ given by
\begin{align}\label{eq:V}
V(\ve{x}):=\tfrac12\sum_{\{i,j\}\in\E}a_{ij}\|\ve[j]{x}-\ve[i]{x}\|^2_2,
\end{align}
where $a_{ij}=a_{ji}\in[0,\infty)$.
Clearly, $V(\ve{x})=0$ if and only if $\ve{x}=(\ve[i]{x})_{i=1}^N\in\mathcal{C}$, \ie if there is no disagreement

The CSS that we study in this paper, Algortihm \ref{algo:my}, is the gradient descent flow of \eqref{eq:V} on $\M$. The general formulation of this algorithm first appears in \citet{sarlette2009consensus}, although their work is limited to the case when the norm of the states are constant, $\|\ve[i]{x}\|=k$, \ie the case when $\M$ is a subset of a sphere. This includes several important cases such as the Stiefel manifold and Grassmannian manifold. 
\begin{myalgorithm}\label{algo:my}
Setting $\vd[i]{x}=-\mathbf{P}_{\ve[i]{x}}(\ve[i]{u})$ with $\ve[i]{u}=\nabla_i V$ yields a CSS on $\M$ given by
\begin{align}
\vd[i]{x}=-\mathbf{P}_{\ve[i]{x}}(\nabla_iV)&=\mathbf{P}_{\ve[i]{x}}\bigl(\sum_{j\in\Ni} a_{ij}(\ve[j]{x}-\ve[i]{x})\bigr),\label{eq:P}
\end{align}
where $a_{ij}=a_{ji}\in[0,\infty)$ and $\mathcal{N}_i:=\{j\in\mathcal{V}\,|\,\{i,j\}\in\mathcal{E}\}$.
\end{myalgorithm}

Suppose that $\mathbf{P}_{\ve[i]{x}}$ is $C^2$ for all $\ve[i]{x}\in\M$, then the flow of \eqref{eq:P} exists and is unique for all times $t\in\R$ by Proposition \ref{prop:picard}. Nagumo's theorem implies that the system stays on $\M$ for all future time \citep{blanchini2008set}.


For the special case that $\M$ is a hypersurface  there is an explicit expression for $\mathbf{P}_{\ve[i]{x}}$. Using this expression, the dynamics \eqref{eq:P} simplify as
\begin{align}\label{eq:hyper}
\vd[i]{x}&=-\mathbf{P}_{\ve[i]{x}}(\nabla_iV)=\left(\ma{I}-\ve[i]{n}\vet[i]{n}\right)\sum_{j\in\Ni}a_{ij}( \ve[j]{x}-\ve[i]{x}),
\end{align}
where $\ve[i]{n}=\nabla_i c(\ve[i]{x})/\|\nabla_i c(\ve[i]{x})\|$. 

\section{Main results}

\noindent The results of this section are divided into two parts. First, we give an example to show that a CSS on the form \eqref{eq:hyper} can diverge to infinity as suggested by Proposition \ref{prop:diverge}. Moreover we provide a result, Theorem \ref{prop:unbounded}, which contains a condition that excludes this possibility. Second, we show by example that a CSS on a compact manifold can be initialized arbitrarily close to $\mathcal{C}$ but still not converge to it ($\mathcal{C}$ is stable but not asympotically stable). Then we give a result, Theorem \ref{prop:local}, with a condition that excludes this possibility.

\subsection{Divergence on an unbounded hypersurface}



A hypersurface of revolution $\mathcal{H}^2\subset\R^3$ is obtained by rotating a curve around an axis. It has azimuthal (cylindrical) symmetry around the axis, which we assume to be the $z$-axis. The surface $\mathcal{H}^2$ can then be parametrized as
\begin{align}\label{eq:revolution}
\ve{x}:\R^2\rightarrow\R^3:(u,v)\mapsto\begin{bmatrix}
\phi(u)\cos v\\
\phi(u)\sin v\\
\psi(u)
\end{bmatrix},
\end{align}
where $\phi,\psi:\R\rightarrow\R$ are two continuous functions. 

\begin{definition}
For surfaces of revolution in $\R^3$, defined by two given functions $\phi,\psi$, and cycle networks $\mathcal{C}_N$ given by \eqref{eq:cycle}, we define $q$-twisted configurations as
\begin{align*}
\mathcal{T}_q=\Bigl\{(\ve[i]{x})_{i=1}^N\in\M^N\,|\,\ve[i]{x}&=\begin{bmatrix}
\phi(u)\cos (\theta+2\pi qi/N)\\
\phi(u)\sin (\theta+2\pi qi/N)\\
\psi(u)
\end{bmatrix},\\
\theta&\in\R,\,u\in\R,\,\forall\,i\in\V\Bigr\}
\end{align*}
\end{definition}



The next example shows that the gradient descent flow diverges on some manifolds, even as it converges to consensus, \ie it is possible that $\lim_{t\rightarrow\infty}\ve[i]{x}(t)-\ve[j]{x}(t)=\ve{0}$ for all $\{i,j\}\in\E$ while $\lim_{t\rightarrow\infty}\|\ve[i]{x}(t)\|=\infty$ as suggested by Proposition \ref{prop:diverge}. Moreover, this example has the property that both $d((\ve[i]{x})_{i=1}^N,\mathcal{C})\rightarrow0$ and $d((\ve[i]{x})_{i=1}^N,\mathcal{T}_1)\rightarrow0$ as $n\rightarrow\infty$.

\begin{example}
\label{ex:pseudosphere}
Let $\mathcal{H}^2$ be the pseudosphere, which is a hypersurface of constant negative curvature in $\R^3$, see Fig \ref{fig:pseudosphere}. A parametrization of the pseudosphere is given by
\begin{align}\label{eq:traxtic}
\ve{y}:\R^2\rightarrow\M\subset\R^3:(u,v)\mapsto\begin{bmatrix}
\sech u\cos v\\
\sech u\sin v\\
u-\tanh u
\end{bmatrix},
\end{align}
where $v\in(-\pi,\pi]$ and $u\in\R$. For $u\in(0,\infty)$ an outward pointing unit normal of the pseudosphere is given by
\begin{align*}
\ve{n}=\begin{bmatrix}
\tanh u\cos v\\
\tanh u\sin v\\
\sech u
\end{bmatrix}.
\end{align*}

\begin{figure}[htb]
	\centering
	\includegraphics[trim=4cm 2cm 4cm 2cm, clip,width=0.23\textwidth]{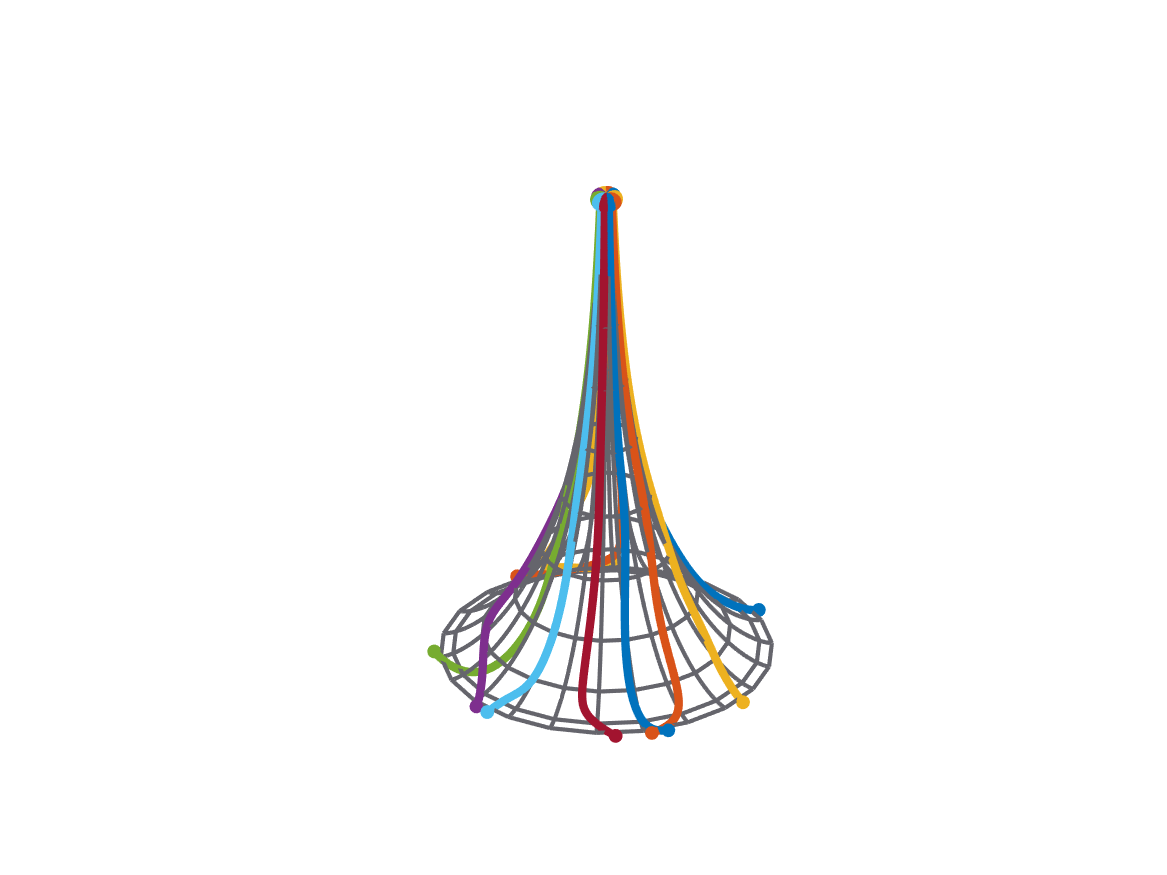}
	\caption{\label{fig:pseudosphere}A CSS of 10 agents on the pseudosphere diverges to infinity.} 
\end{figure}

Since the pseudosphere is a surface of revolution around the $z$-axis, we can distribute the agents equidistantly at a constant $u$ value to obtain a $1$-twisted configuration:
\begin{align*}
\ve[i]{x}(u)&=\begin{bmatrix}
\sech u\cos 2\pi i/N\\
\sech u\sin 2\pi i/N\\
u-\tanh u
\end{bmatrix},
\end{align*}
$(\ve[i]{x})_{i=1}^N\in\mathcal{T}_1$. The agent are spread out on a circle of radius $\sech u$ (all agents have the same $u$ value) and are connected by the cycle graph \eqref{eq:cycle}.

The input $\ve[i]{u}$ satisfies
\begin{align*}
\ve[i]{u}=&{}\begin{bmatrix}
\sech u\cos 2\pi (i-1)/N\\
\sech u\sin 2\pi (i-1)/N\\
u-\tanh u
\end{bmatrix}+\begin{bmatrix}
\sech u\cos 2\pi (i+1)/N\\
\sech u\sin 2\pi (i+1)/N\\
u-\tanh u
\end{bmatrix}\\
&-2\begin{bmatrix}
\sech u\cos 2\pi i/N\\
\sech u\sin 2\pi i/N\\
u-\tanh u
\end{bmatrix}\\
={}&2(\cos2\pi/N-1)\begin{bmatrix}
\sech u\cos 2\pi i/N\\
\sech u\sin 2\pi i/N\\
0
\end{bmatrix}\\
={}&\frac{2(\cos2\pi/N-1)}{\sinh u}\ve[i]{n}-\frac{2(\cos2\pi/N-1)\sech u}{\sinh u}\ve[3]{e},
\end{align*}
where $\ve[i]{n}$ denotes the normal of the pseudosphere at $\ve[i]{x}$ and $\ve[3]{e}=[0 \,0 \,1]\mtr$. It follows that 
\begin{align}
\vd[i]{x}&=(\ma[3]{I}-\ve[i]{n}\vet[i]{n})\ve[i]{u}=(\ma[3]{I}-\ve[i]{n}\vet[i]{n})f(u)\ve[3]{e},\label{eq:reduced}\\
f(u)&=\frac{2(1-\cos2\pi/N)}{\cosh u\sinh u}=\frac{4(1-\cos2\pi/N)}{\sinh(2u)}.\nonumber
\end{align}

Introduce a Lyapunov function $U=\vet[3]{e}\ve[i]{x}=u-\tanh(u)$ with derivatives
\begin{align*}
\dot{U}&=f(u)(1-(\vet[3]{e}\ve[i]{n})^2)=f(u)\tanh^2(u),\\
\dot{U}&=\dot{u}-\dot{u}(1-\tanh^2(u))=\dot{u}\tanh^2(u).
\end{align*}
We thus find $\dot{u}=f(u)$. This equation can be solved for
\begin{align*}
u(t)&=\cosh^{-1}(8(1-\cos 2\pi/N)t+\cosh(u(0))).
\end{align*}
Note that $u(t)\rightarrow\infty$ as $t\rightarrow\infty$, see also Fig. \ref{fig:pseudosphere}. \hfill$\blacklozenge$
\end{example}




\subsection{Boundedness on boundaries of convex sets}





The following result states conditions under which divergence like that we saw in Example \ref{ex:pseudosphere} is not possible.

\begin{theorem}\label{prop:unbounded}
If a hypersurface $\mathcal{H}^n$ is the boundary $\partial\mathcal{K}$ of a convex set $\mathcal{K}\subset\R^{n+1}$, then the CSS \eqref{eq:hyper} on $\M$ satisfies
\begin{align*}
\lim_{t\rightarrow\infty}\|\ve[i]{x}(t)\|<\infty
\end{align*}
for all initial conditions on $(\mathcal{H}^n)^N$ and all $i\in\V$.
\end{theorem}

\begin{pf} 
First we establish a property of unbounded sequences in convex sets. Translate $\mathcal{K}$ so that $\ve{0}\in\mathcal{K}$. Let $\ve{s}:=(\ve[i]{s})_{i=1}^\infty$	be any sequenence such that  $\ve[i]{s}\in\mathcal{K}$ and    $\lim_{i\rightarrow\infty}\|\ve[i]{s}\|=\infty$. Take a subsequence $\ve{s}^\prime:=(\ve[i]{s}^\prime)_{i=1}^\infty$ of $\ve{s}$ such that $\|\ve[i]{s}^\prime\|\geq i$. Form a sequence $\ve{r}:=(\ve[i]{r})_{i=1}^\infty$ on $\mathcal{S}^{n}$ by setting $\ve[i]{r}:=\ve[i]{s}^\prime/\|\ve[i]{s}^\prime\|$. Since $\ve{r}$ is a sequence on a bounded set, by the Bolzano-Weierstrass theorem, there exists a $\ve{d}\in\mathcal{S}^n$ such that $\ve{r}$ has a convergent subsequence $\ve{r}^\prime$ with $\lim_{i\rightarrow\infty}\ve[i]{r}^\prime=\ve{d}$. Consider any point $\alpha\ve{d}$ for some fixed $\alpha\in(0,\infty)$. Because $\ve{0},\ve[i]{s}^\prime\in\mathcal{K}$ and $\mathcal{K}$ is convex, for $i\geq\alpha$ we have $\alpha/\|\ve[i]{s}^\prime\|\in[0,1]$ whereby
\begin{align*}
\alpha\ve[i]{s}^\prime/\|\ve[i]{s}^\prime\|=\alpha/\|\ve[i]{s}^\prime\|\,\ve[i]{s}^\prime+(1-\alpha/\|\ve[i]{s}^\prime\|)\,\ve{0}\in\mathcal{K}.
\end{align*}
Moreover, since $\mathcal{K}$ is closed, $\alpha\ve{d}=\lim_{i\rightarrow\infty}\alpha\ve[i]{s}^\prime/\|\ve[i]{s}^\prime\|\in\mathcal{K}$. It follows that $\{\alpha\ve{d}\,|\,\alpha\in[0,\infty)\}\subset\mathcal{K}$.

Suppose that $\|\ve[i]{x}(t)\|\rightarrow\infty$ as $t\rightarrow\infty$. By the reasoning of the previous paragraph, there exists a $\ve{d}\in\mathcal{S}^n$ such that $\sup_{i\in\N}\langle\ve[i]{x}(t),\ve{d}\rangle=\infty$. Let $\mathcal{P}_{\alpha}$ be a family of parallel affine planes such that $\alpha \ve{d}\in\mathcal{P}_{\alpha}$ and $\ve{d}\perp\ve{p}-\ve{q}$ for all $\ve{p},\ve{q}\in\mathcal{P}_{\alpha}$, \ie $\ve{d}$ is a normal of $\mathcal{P}_\alpha$. Note that $\mathcal{P}_\alpha$ separates $\mathcal{K}$ into two subsets $\mathcal{K}_0$ and $\mathcal{K}_\infty$:
\begin{align*}
\ve{0}\in\mathcal{K}_0,\quad\mathcal{K}_0\cup\mathcal{K}_\infty=\mathcal{K},\quad \mathcal{K}_0\cap\mathcal{K}_\infty\subset\mathcal{P}_\alpha,
\end{align*}
and $\mathcal{K}_\infty$ is unbounded. Note that $\|\ve{k}\|\geq\alpha$ for all $\ve{k}\in\mathcal{K}_\infty$. Since $\|\ve[i]{x}(0)\|$ is finite we can pick an $\alpha\in(0,\infty)$ such that $\{\ve[1]{x}(0),\ldots,\ve[N]{x}(0)\}\subset\mathcal{K}_{0}$. 

Note that for $\ve[i]{x}$ to diverge along the direction of $\ve{d}$, it must first pass through $\mathcal{P}_\alpha$. Let $i$ be the index of the agent whose state $\ve[i]{x}$  passes through $\mathcal{P}_\alpha$ first, \ie at the earliest time $\tau$. At that time $\tau$, we note that
\begin{align}
\langle\ve{d},\vd[i]{x}\rangle={}&\langle\ve{d},(\ma{I}-\ve[i]{n}\vet[i]{n})\sum_{j\in\Ni}a_{ij}(\ve[j]{x}-\ve[i]{x})\rangle\nonumber\\
={}&\langle\ve{d},\sum_{j\in\Ni}a_{ij}(\ve[j]{x}-\ve[i]{x})\rangle-\nonumber\\
&\langle\ve[i]{n},\ve{d}\rangle\langle\ve[i]{n},\sum_{j\in\Ni}a_{ij}(\ve[j]{x}-\ve[i]{x})\rangle.\label{eq:ineq}
\end{align}

We claim that the last expression in \eqref{eq:ineq} is negative. To see this, note that $\langle\ve{d},\sum_{j\in\Ni}a_{ij}(\ve[j]{x}-\ve[i]{x})\rangle\leq 0$ since $\ve{d}$ is normal to $\mathcal{P}_\alpha$ and points towards $\mathcal{K}_\infty$ whereas $\ve[j]{x}-\ve[i]{x}$ points away from $\mathcal{K}_\infty$ due to $\ve[j]{x}$ belonging to $\mathcal{K}_0$. The inequality follows from $\ve{d}$ and $\ve[j]{x}-\ve[i]{x}$ belonging to different half-spaces. By the Hahn-Banach separation theorem, since $\mathcal{K}$ is convex, the tangent plane $\mathsf{T}_{\ve[i]{x}}\mathcal{K}$ separates $\R^{n+1}$ into two halfspaces, one that contains $\mathcal{K}$ and one that does not. Moreover, $\ve[i]{n}$ is normal to $\smash{\mathsf{T}_{\ve[i]{x}}}\mathcal{K}$ and we may assume that it points towards the halfspace which does not contain $\mathcal{K}$. As such $\langle\ve[i]{n},\ve{k}-\ve[i]{x}\rangle\leq0$ for all $\ve{k}\in\mathcal{K}$, including $\ve[j]{x}\in\mathcal{K}$. This also yields $\langle\ve[i]{n},\ve{d}\rangle\leq0$ since $\ve{d}=\ve{k}-\ve[i]{x}$ for $\ve{k}=\ve[i]{x}+\ve{d}\in\mathcal{K}$.


In summary, $\ve[i]{x}$ cannot pass through $\mathcal{P}_\alpha$ since $\langle\ve{d},\ve[i]{x}\rangle$ is decreasing, contradicting our assumption about $\tau$.\hfill$\blacksquare$\end{pf}

%
%



\subsection{Divergence on a compact hypersurface}


Our next example shows that it is possible to initialize a CSS on a hypersurface arbitrarily close to $\mathcal{C}$ with having $\lim_{t\rightarrow\infty}d(\ve{x}(t),\mathcal{C})=0$, \ie $\mathcal{C}$ is not asymptotically stable.

\begin{example}
Construct a telescope like hypersurface by glueing together infinitely many cylinders. Denote the cylinders by $k\in\N$. The length and diameter of the $k$th cylinder is $2^{-k}$. Note that each pair of subsequent cylinders $k$, $k+1$ can be glued together in an arbitrarily smooth manner (\ie $C^\infty$) using bump functions. However, they cannot be glued together using analytic functions. The final result will hence not be an analytic manifold, as is important in Theorem \ref{prop:local}.

We derive the dynamics of a CSS on the telescope. First, we find an implicit characterization of the $k$th cylinder
\begin{align*}
\left\{\begin{bmatrix}
x\\
y\\
z\\
\end{bmatrix}\in\R^3\,|\,x^2+y^2=2^{-k},\,\sum_{j=1}^{k-1}2^{-j}\leq z\leq \sum_{j=1}^{k}2^{-j}\right\}.
\end{align*}
Moreover, the normal of the $k$th cylinder is given by
\begin{align*}
\ve{n}=\frac{1}{\sqrt{x^2+y^2}}
\begin{bmatrix}
x\\y\\0
\end{bmatrix}
=2^{k/2}\begin{bmatrix}
x\\y\\0
\end{bmatrix}.
\end{align*}
Assume the agents are connected by the cycle graph \eqref{eq:cycle} and set $a_{ij}=1$. It follows that the dynamics \eqref{eq:hyper} are
\begin{align*}
\vd[i]{x}&=\Bigl(\ma[3]{I}-2^{k}\begin{bmatrix}
x_i\\
y_i\\
0
\end{bmatrix}\begin{bmatrix}
x_i & y_i & 0
\end{bmatrix}\Bigr)(\ve[i-1]{x}+\ve[i+1]{x}-2\ve[i]{x}),
\end{align*}
which decouples into
\begin{align*}
\begin{bmatrix}
\dot{x}_i\\
\dot{y}_i
\end{bmatrix}={}&\Bigl(\ma[3]{I}-2^{k}\begin{bmatrix}
x_i\\
y_i
\end{bmatrix}\begin{bmatrix}
x_i & y_i
\end{bmatrix}\Bigr)\cdot\\
&\Bigl(\begin{bmatrix}
x_{i-1}\\
y_{i-1}
\end{bmatrix}+\begin{bmatrix}
x_{i+1}\\
y_{i+1}
\end{bmatrix}-2\begin{bmatrix}
x_{i}\\
y_{i}
\end{bmatrix}\Bigr)\\
\dot{z}_i={}&z_{i+1}+z_{i-1}-2z_i.
\end{align*}

Since the $z$-dynamics are a linear CSS on $\R$, it follows that $\sum_{j=1}^{k-1}2^{-j}\leq z_i\leq\sum_{j=1}^{k}2^{-j}$ for all $i$ is an invariant set. This means the agents will stay on a single cylinder if they are initialized on it. The dynamics of $x_i$ and $y_i$ are the CSS \eqref{eq:hyper} on the circle $\mathcal{S}^1=\{\ve{y}\in\R^2\,|\,\|\ve{y}\|^2=1\}$ with the cycle graph \eqref{eq:cycle}. This dynamics is the well-known homogeneous Kuramoto model on the circle in Cartesian coordinates. In particular, it holds that $1$-twisted equilibria are asymptotically stable for this CSS \citep{wiley2006size}.

The point of this example is to show that if we initialize the agents in a $1$-twisted equilibrium on the $k$th cylinder with $a_{ij}=1$, then the value of the disagrement function is
\begin{align*}
V(\ve{x})&=\tfrac12\sum_{i,j\in\E}a_{ij}\|\ve[j]{x}-\ve[i]{x}\|^2_2=2^{-k+1}N\sin\pi/N,
\end{align*}
where we used the fact that the agents are equidistantly spread over a circle and the chord length formula. It is now clear that for any $\varepsilon>0$ we can position the agents on a cylinder with a $k$ sufficiently large that the disagrement function satisfies $V(\ve{x})<\varepsilon$. In particular, the agents can start out arbitrarily near $\mathcal{C}$, \ie near $V(\ve{x})=0$, without converging to $\mathcal{C}$ since they converge to $\mathcal{T}_1$ instead.\hfill$\blacklozenge$
\end{example}

Another finding of note is that there is a hypersurface on which the state $\ve{x}$ of the CSS approaches a limit set rather than a singleton. This result can be shown by building a CSS out of the so-called Mexican hat function \citep{absil2005convergence}. It is omitted due to space constraints. 

\subsection{Convergence on closed analytic manifolds}

\noindent Finally, we give a condition such that the CSS \eqref{eq:P} converges to synchronization if $\ve{x}$ is initialized sufficiently close to $\mathcal{C}$, \ie conditions such that $\mathcal{C}$ is asymptotically stable.

\begin{theorem}\label{prop:local}
Let $\M\subset\R^m$ be a closed analytic Riemannian manifold. Consider the gradient descent flow $\vd{x}=-\mathbf{P}_{\ve{x}}\nabla V$,
where $\ve{x}=(\ve[i]{x})_{i=1}^N$ and
\begin{align*}
V(\ve{x})&:=\tfrac12\sum_{\{i,j\}\in\E}\|\ve[j]{x}-\ve[i]{x}\|^2.
\end{align*}
\ie the system given by Algorithm \ref{algo:my}. For each initial condition, the system converges to a singleton set. Moreover, the synchronization manifold $\mathcal{C}:=\{(\ve[i]{x})_{i=1}^N\in\M^N\,|\,\ve[i]{x}=\ve[j]{x},\,\forall\,\{i,j\}\in\mathcal{E}\}$ is asymptotically stable.
\end{theorem}

\begin{pf}
Since $\M$ is closed, the system cannot diverge to infinity. Convergence to a singleton set follows from a result on analytic gradient descent flow \citep{lageman2007convergence}.

The potential function of a gradient descent flow decreases with time, $\dot{V}(\ve{x})=\langle\nabla V(\ve{x}),\vd{x}\rangle=-\|\nabla V(\ve{x})\|^2$. Since $V(\ve{x})\geq 0$ with $V(\ve{x})=0$ if and only if $\ve{x}\in\mathcal{C}$, we can take $V(\ve{x})$ as a Lyapunov function and conclude that $\mathcal{C}$ is stable. 

Since $\M$ is closed, the gradient descent flow converges to a connected component of the set of critical points of $V$ \citep{helmke2012optimization}. Any sublevel set of $V$ is forward invariant. Moreover, all sublevel sets contain $\mathcal{C}=\{\ve{x}\in\M^N\,|\,V(\ve{x})=0\}$. Let $\mathcal{Q}$ denote the set of equilibria of \eqref{eq:P} that do not belong to $\mathcal{C}$. If there is an open sublevel set of $V$ which does not intersect $\mathcal{Q}$, then there is an open neighborhood of $\mathcal{C}$ from which $\ve{x}$ converges to $\mathcal{C}$. 

Since $V(\ve{x})$ is analytic it satisfies the \L{}ojasiewicz inequality on Riemannian manifolds \citep{kurdyka2000proof}. For every $\ve{x}\in\mathcal{C}$ there is an open ball $\mathcal{B}(\ve{x})$, an $\alpha<1$, and a $k>0$ such that $V(\ve{y})^\alpha\leq k\|\nabla V(\ve{y})\|$ for all $\ve{y}\in\mathcal{B}(\ve{x})$. If $\ve{y}\in\mathcal{Q}$, then $\nabla V(\ve{y})=\ve{0}$ whereby $V(\ve{y})=0$. However, this implies  $\ve{y}\in\mathcal{C}$, a contradiction. Hence $\mathcal{Q}\cap\mathcal{B}(\ve{x})=\emptyset$.

Consider the value of $q=\inf_{\ve{x}\in\mathcal{Q}}V(\ve{x})$. If $q=0$, then there is a sequence $\{\ve[k]{x}\}_{k=1}^\infty$ such that $\lim_{k\rightarrow\infty}V(\ve[k]{x})=0$. Since $\M$ is a closed manifold embedded in $\R^n$, by the Bolzano-Weierstrass theorem, the sequence $\{\ve[k]{x}\}_{k=1}^\infty$ has a subsequence which converges to some $\ve{y}\in\M$. Moreover, $V(\ve{y})=0$ whereby $\ve{y}\in\mathcal{C}$. For each $\varepsilon>0$ there must be a $\ve{z}(\varepsilon)\in\mathcal{Q}$ (an element of the convergent subsequence) such that $\|\ve{y}-\ve{z}(\varepsilon)\|<\varepsilon$. This contradicts  $\mathcal{Q}\cap\mathcal{B}(\ve{y})=\emptyset$. Hence $q>0$ and all trajectories that start in the level set $\{\ve{x}\in\M\,|\,V(\ve{x})<q\}$ converges to $\mathcal{C}$. \hfill$\blacksquare$\end{pf}

\section{Conclusions}

The results of this paper are given by examples which show how consensus seeking systems (CSS) can go wrong, figuratively speaking, and theorems that provide the assumptions needed to overcome these deficiencies. The results are based on a formulation of CSS as gradient descent flows of a disagreement function. This setting is chosen since it  gives us access to multiple powerful tools including the \L{}ojasiewicz gradient inequality \citep{kurdyka2000proof}. A limitation of this CSS framework is that it is restricted to systems with undirected interaction graphs, owing to a symmetry property of gradients. Another approach based on positivity and monotone system design provides convergence results for directed graphs on Lie groups \citep{mostajeran2018positivity}. This builds on earlier work \citep{sepulchre2010consensus} that utilizes a measure of the diameter of the convex hull of agents states as a Lyapunov function \citep{blondel2005convergence}.

\bibliographystyle{plain} 
\bibliography{autosam}






\end{document}

%% file: myIncludes.tex
\usepackage{graphicx}

\usepackage{subdepth} 

\usepackage{graphics} 
\usepackage{epsfig} 

\usepackage{amssymb} 
\usepackage{amsmath,psfrag,color,graphicx} 

\usepackage{graphicx}
\usepackage{quoting} 
\usepackage{lettrine} 
\usepackage{tocloft} 
\usepackage{titletoc}
\usepackage{fmtcount}
\usepackage{calc} 
\usepackage{cite} 
\usepackage{multicol} 
\usepackage{centernot}
\usepackage{nicefrac}
\usepackage{dsfont}
\usepackage{booktabs}
\usepackage{mathtools, cuted}

\usepackage{lipsum}
\usepackage{psfrag}
\usepackage[export]{adjustbox} 

\usepackage{mathtools}

\PassOptionsToPackage{cmyk}{xcolor}
\usepackage[cmyk]{xcolor}
\selectcolormodel{cmyk}

\usepackage{textcomp}

\usepackage{txfonts}
\usepackage{lmodern}
\usepackage{times}

\DeclareFontEncoding{LS1}{}{}
\DeclareFontSubstitution{LS1}{stix}{m}{n}

\makeatletter
\newcommand*\textmathversion{\csname textmv@\math@version\endcsname}
\newcommand*\textmv@normal{m}
\newcommand*\textmv@bold{b}
\makeatother

\newcommand{\ve}[2][]{\ensuremath{\boldsymbol{\mathrm{#2}}}_{#1}}
\newcommand{\vet}[2][]{\ensuremath{\smash{\boldsymbol{\mathrm{#2}}^{\!\top}_{#1}}}}
\newcommand{\vd}[2][]{\ensuremath{\dot{\boldsymbol{\mathrm{#2}}}_{#1}}}

\newcommand{\ma}[2][]{\ensuremath{\boldsymbol{\mathrm{#2}}}_{#1}}

\DeclareMathOperator{\grad}{\ensuremath{\mathrm{grad}}}

\newcommand{\R}{\ensuremath{\mathds{R}}}

\newcommand{\N}{\ensuremath{\mathds{N}}}

\newcommand{\ie}{\textit{i.e.}, }

\newcommand{\eg}{\textit{e.g.}, }

\newcommand{\mtr}{\hspace{-0.3mm}\ensuremath{^\top}}

\newcommand\raiseT[2]{\setbox0\hbox{$#1{#2}$}\raise\dp0\box0}

\newcommand{\V}{\ensuremath{\mathcal{V}}}

\newcommand{\E}{\ensuremath{\mathcal{E}}}

\newcommand{\M}{\ensuremath{\mathcal{M}}}

\newcommand{\etc}{\emph{etc.\ }}

\DeclareMathOperator{\sech}{\ensuremath{\mathrm{sech}}}

\newcommand{\Ni}{\ensuremath{\mathcal{N}_i}}

\makeatletter
\newcommand{\raisemath}[1]{\mathpalette{\raisem@th{#1}}}
\newcommand{\raisem@th}[3]{\raisebox{#1}{$#2#3$}}
\makeatother

\newcommand{\ts}[2][]{\ensuremath{\mathsf{T}_{#2}#1}}

\definecolor{kthbluergb}{RGB}{25,84,166}
\definecolor{kthbluecmyk}{cmyk}{1,0.55,0,0}

\definecolor{kthblueA}{RGB}{25,84,166}
\definecolor{kthblueB}{RGB}{46,124,192}
\definecolor{kthblueC}{RGB}{112,153,209}
\definecolor{kthblueD}{RGB}{164,186,225}
\definecolor{kthblueE}{RGB}{211,220,241}

\newcounter{counter} 

\newtheorem{theorem}[counter]{Theorem}

\newtheorem{proposition}[counter]{Proposition}

\newtheorem{definition}[counter]{Definition}
\newtheorem{example}[counter]{Example}

\newtheorem{myalgorithm}[counter]{Algorithm}


\newcounter{parentnumber}
